\documentclass[11pt,a4paper]{article}
\usepackage{mathrsfs}
\usepackage{pb-diagram}
\newtheorem{theorem}{Theorem}[section]
\newtheorem{lemma}[theorem]{Lemma}
\newtheorem{corollary}[theorem]{Corollary}
\newtheorem{proposition}[theorem]{Proposition}

\newtheorem{definition}{Definition}[section]
\newtheorem{example}{Example}[section]

 \input amssym.def
\input amssym.tex

\def\bc{\begin{center}}
\def\ec{\end{center}}

\begin{document}
\title{Representations of skew group algebras induced from isomorphically invariant modules over path algebras\thanks{Project supported by the Program for New Century Excellent Talents in University
(No.04-0522) and the National Natural Science Foundation of China
(No.10571153) }}

\author{Mianmian Zhang$^{1,2}\thanks{zhmm1216@yahoo.com.cn}\;\;\;\;\;\;\;$ Fang Li$^1$\thanks{fangli@cms.zju.edu.cn} \\{$^1$Department of Mathematics, Zhejiang University,
Hangzhou, Zhejiang 310027, China}\\{$^2$Department of Mathematics,
Kansas State University, Manhattan, KS 66506, USA}}
 \maketitle

\begin{abstract}
Suppose that $Q$ is a connected quiver without oriented cycles and
$\sigma$ is an automorphism of $Q$. Let $k$ be an algebraically
closed field whose characteristic does not divide the order of the
cyclic group $\langle\sigma\rangle$.

The aim of this paper is to investigate the relationship between
indecomposable $kQ$-modules and indecomposable
$kQ\#k\langle\sigma\rangle$-modules. It has been shown by Hubery
that any $kQ\#k\langle\sigma\rangle$-module is an isomorphically
invariant $kQ$-module, i.e., ii-module (in this paper, we call it
$\langle\sigma\rangle$-equivalent $kQ$-module), and conversely any
$\langle\sigma\rangle$-equivalent $kQ$-module induces a
$kQ\#k\langle\sigma\rangle$-module. In this paper, the authors
prove that a $kQ\#k\langle\sigma\rangle$-module is indecomposable
if and only if it is  an indecomposable
$\langle\sigma\rangle$-equivalent $kQ$-module. Namely, a method is
given in order to induce all indecomposable
$kQ\#k\langle\sigma\rangle$-modules from all indecomposable
$\langle\sigma\rangle$-equivalent $kQ$-modules.  The number of
non-isomorphic  indecomposable $kQ\#k\langle\sigma\rangle$-modules
induced from the same indecomposable
$\langle\sigma\rangle$-equivalent $kQ$-module is given. In
particular, the authors give the relationship between
indecomposable $kQ\#k\langle\sigma\rangle$-modules and
indecomposable $kQ$-modules in the cases of indecomposable simple,
projective and injective modules.

\end{abstract}
\section{Introduction}

There is a lot of literature on smash product algebras $A\#H$ and
crossed product algebras $A\#_{\sigma}H$, and on their
relationships with the algebra $A^{H}$, whose elements are those
elements of $A$ left fixed by $H$. Much work has been done to
determine which properties of $A$ are inherited by
$A\#_{\sigma}H$. The works about the relationships among these
algebras were motivated by the development
 of the Galois theory of noncommutative algebras.

It is important to totally understand the relationships among
representations of $A$, $A\#H$ and $A\#_{\sigma}H$. It has been
proven in \cite{L,LZ} that the representation types of $A$, $A\#H$
and $A\#_{\sigma}H$ are the same if $A$ is a finite dimensional
algebra and $H$ is a finite dimensional semisimple and
cosemisimple Hopf algebra over an algebraically closed field.
These results were used to classify finite dimensional basic Hopf
algebras through their representation types in \cite{L,Liu, Liu
2,LL}. In particular, each $A\#H$-module (or respectively,
$A\#_{\sigma}H$-module) is an $A$-module, but not all $A$-modules
induce $A\#H$-modules (or respectively, $A\#_{\sigma}H$-modules).
Thus, some questions arise, such as the following examples:
\begin{itemize}
 \item[$(i)$] What
kind of $A$-modules can induce $A\#H$-modules (or respectively,
$A\#_{\sigma}H$-modules)?
 \item[$(ii)$] If an $A$-module can induce
$A\#H$-modules (or respectively, $A\#_{\sigma}H$-modules), how
many non-isomorphic classes of such induced $A\#H$-modules (or
respectively, $A\#_{\sigma}H$-modules) exist?
\end{itemize}

In \cite{L,LZ}, we have proven that for finite dimensional algebra
$A$, finite dimensional Hopf algebra $H$ such that $H$ and its
dual $H^{*}$ are both semisimple, then  $A$, $A\#H$ and
$A\#_{\sigma}H$ have the same representation type. This result
allows us the possibility to discuss these questions.

Hubery in \cite{H,H2} constructed  the dual quiver with
automorphism ($\widetilde{Q},\tilde{\sigma}$), where
$\widetilde{Q}$ is the Ext-quiver of $kQ\#k\langle\sigma\rangle$
and $\tilde{\sigma}$ is the automorphism of $k\tilde{Q}$ induced
from an admissible automorphism $\sigma$. Here, the  admissible
automorphism $\sigma$ means that $Q$ has no arrow connecting two
vertices in the same $\sigma$-orbit, and $k$ is an algebraically
closed field of characteristic not dividing the order of
$\langle\sigma\rangle$.
 Hubery used the dual quiver $(\widetilde{Q},\tilde{\sigma})$ to
prove the generalization of Kac's Theorem. During the
construction, Hubery defined the isomorphically invariant module,
i.e. ii-module (in this paper, we call it
$\langle\sigma\rangle$-equivalent $kQ$-module), and proved
 that any $kQ\#k\langle\sigma\rangle$-module is an
 ii-module and conversely any
ii-module induces a $kQ\#k\langle\sigma\rangle$-module. These
works by Hubery are crucial for us to answer the above questions
in the special case when $A$ is a path algebra $kQ$, and
$H=k\langle\sigma\rangle$, a cyclic group algebra with $\sigma\in
Aut(Q)$.

In this paper, we investigate the relationship between
indecomposable modules over the path algebra $kQ$ and the skew
group algebra $kQ\#k\langle\sigma\rangle$ respectively, where $k$
is an algebraically closed field with the characteristic not
dividing the order of $\sigma$, $Q$ is connected and without
oriented cycles, and $\sigma\in Aut(Q)$. We prove that a
$kQ\#k\langle\sigma\rangle$-module is indecomposable if and only
if it is  an indecomposable $\langle\sigma\rangle$-equivalent
$kQ$-module. Namely, a method is given in order to induce all
indecomposable $kQ\#k\langle\sigma\rangle$-modules from each
indecomposable $\langle\sigma\rangle$-equivalent $kQ$-module. The
number of non-isomorphic indecomposable
$kQ\#k\langle\sigma\rangle$-modules induced from the same
indecomposable $\langle\sigma\rangle$-equivalent $kQ$-module is
given.

In this paper, assume that all the modules are unital and finitely
generated and that $k$ is always an algebraically closed field.
All the concepts and notations on  Hopf algebra and crossed
product algebra can be found in \cite{M}. We fix the notation
$\Delta(h)=\sum_{(h)}h_{1}\otimes h_{2}$ under the Sweedler
meaning. If $H$ is a group algebra $kG$, $\Delta(g)=g\otimes g,
\varepsilon(g)=1,$ for all $g\in G$. In particular, we give the
definition of smash product as follows:

\begin{definition} Let $H$ be a Hopf algebra. An algebra $A$ is a {\em (left) $H$-module
algebra} if for all $h\in H, a,b\in A$,

\begin{itemize}

 \item[${(1)}$] A is a (left) $H$-module, via
$h\otimes a\mapsto h\cdot a,$

\item[$(2)$]$h\cdot (ab)=\sum_{(h)} (h_{1}\cdot a)(h_{2}\cdot b)$,
\item[$(3)$] $h\cdot 1_{A}=\varepsilon (h)1_{A}.$
\end{itemize}

\end{definition}

\begin{definition} Let $A$ be a left $H$-module algebra. The
{\em smash product algebra} $A\#H$ is defined by satisfying that
\begin{itemize}
\item[$(1)$] as  a $k$-space, $A\#H=A\otimes H$;
 \item[$(2)$] the multiplication is
given by $(a\#h)(b\#l)=\sum_{(h)}a(h_{1}\cdot b)\#h_{2}l$ for all
$a,b\in A, h,l\in H$.

\end{itemize}

 We write that by $a\#h$ the element $a\otimes h\in A\#H$.
\end{definition}

\section{Isomorphically Invariant $kQ$-modules}

Suppose that $Q=(Q_0,Q_1)$ is a {\em quiver} given by the
\emph{vertex set} $Q_{0}$ and the \emph{arrow set } $Q_{1}$ . For
an arrow $\alpha \in Q_{1}$, the vertex $s(\alpha)$ is the
\emph{start vertex of} $\alpha$ and the vertex $t(\alpha)$ is the
\emph{end vertex of} $\alpha$, and we draw
$s(\alpha)\stackrel{\alpha}\rightarrow t(\alpha)$. A \emph{path }
in $Q$ is ($b|\alpha_{n}\cdots\alpha_{1}|a$), where $\alpha_{i}\in
Q_{1}$, for $i=1,\cdots,n$, and $s(\alpha_{1})=a,
t(\alpha_{i})=s(\alpha_{i+1})$, for $i=1,\cdots,n-1$, and
$t(\alpha _{n})=b$.
 The \emph{length of a path } is the number of arrows in it. To
 each arrow $\alpha$ we can assign an edge $\overline{\alpha}$
 where the orientation is forgotten. A \emph{walk} between two
 vertices $a$ and $b$ is given by
 $(b|\overline{\alpha_{n}}\cdots\overline{\alpha_{1}}|a)$,
 where $a \in \{s(\alpha_{1}),t(\alpha_{1})\},
 b \in
 \{s(\alpha_{n}),t(\alpha_{n})\},$
 and for each $i=1,\cdots,n-1,$
 $\{s(\alpha_{i}),t(\alpha_{i})\}\cap\{s(\alpha_{i+1}),t(\alpha_{i+1})\}\neq\varnothing$.
 A quiver is said to be \emph{connected} if for each pair of
 vertices \emph{a} and \emph{b}, there exists a walk between
 them.

Denote by Rep$Q$ the category  of representations of the quiver
$Q$ over $k$.

It is well-known that a \emph{representation} $X=(X_{i},i\in
Q_{0};\;X_{\rho}: X_{s(\rho)}\rightarrow X_{t(\rho)},\rho \in
Q_{1})$ of $Q$ is given by finite dimensional $k$-vector spaces
$X_{i}$ for all $i\in Q_{0}$ and $k$-linear maps $X_{\rho}:
X_{s(\rho)}\rightarrow X_{t(\rho)}$ for all arrows $\rho \in
Q_{1}$; a \emph{morphism} $\theta: X\rightarrow X^{'}$ is given by
 $k$-linear maps $\theta_{i}: X_{i}\rightarrow X_{i}^{'}$ for $i
\in Q_{0}$ satisfying
$X_{\rho}^{'}\theta_{s(\rho)}=\theta_{t(\rho)}X_{\rho}$ for all
$\rho \in Q_{1}$. The composition of $\theta$ with another
morphism $\phi: X^{'}\rightarrow X^{''}$ is defined by
$(\phi\theta)_{i}=\phi_{i}\theta_{i}$ for all $i \in Q_{0}$.

Let $Q$ be a finite quiver, that is, $\mid Q_0\mid$ and $\mid
Q_1\mid$ are both finite. Denote by mod$kQ$ the category of finite
generated $kQ$-modules. All through the paper, $Q$ is a connected
quiver without oriented cycles.

 For a $kQ$-module $\mathcal{X}$, define a
representation $X$ with the $k$-vector spaces
$X_{i}=e_{i}\mathcal{X}$ for all vertices $i\in Q_0$ and the
linear maps $X_{\rho}$ for all arrows $\rho\in Q_1$ satisfying
$X_{\rho}(x)=\rho x=e_{t(\rho)}\rho x \in X_{t(\rho)}$  for $x\in
X_{s(\rho)}$.
 Conversely, for a representation $X$ of $Q$, define a $kQ$-module $\mathcal{X}$ via
 $\mathcal{X}=\bigoplus_{i\in Q_0}X_{i}$ with actions of paths $\rho_{1}\cdots\rho_{m}$ satisfying
  $\rho_{1}\cdots\rho_{m}x=\varepsilon_{t(\rho_{1})}X_{\rho_{1}}\cdots
X_{\rho_{m}}\pi_{s(\rho_{m})}(x)$ and
$e_{i}x=\varepsilon_{i}\pi_{i}(x)$ for
 the canonical
maps $X_{i}\stackrel{\varepsilon_{i}}\rightarrow
\mathcal{X}\stackrel{\pi_{i}}\rightarrow X_{i}$. Then, as we have
well-known, this correspondence gives a pair of mutually
quasi-invertible functors between Rep$Q$ and mod$kQ$, that is,

\begin{theorem}$(\cite{ASS,ARS})$ \label{thm 2.1}  For  a finite  quiver  $Q$ over a field $k$, the categories Rep$Q$ and mod$kQ$ are equivalent.
\end{theorem}
The correspondence, given above between objects of Rep$Q$ and
mod$kQ$, will be useful for our further discussion.

From now on, we let $Q$ be a connected quiver and without oriented
cycles, $\sigma\in Aut(Q)$ and $k$ is an algebraically closed
field with the characteristic of not dividing the order of
$\sigma.$

It is easy to extend  $\sigma$ linearly to the whole $k$-linear
space $kQ$ as a $k$-automorphism, i.e., $\sigma\in Aut_kkQ$.

Let $\mathcal{X}$ be a $kQ$-module. We define a $kQ$-module
$^{\sigma}\mathcal{X}$ by taking the same underlying vector space
as $\mathcal{X}$ but with the new action:
\[p\cdot x:= \sigma^{-1}(p)x\;\;\; for\;\; p\in kQ.\]
Let $\phi: \mathcal{X}\rightarrow \mathcal{Y}$ be a module
homomorphism, and set $^{\sigma}\phi=\phi$ as a linear map. Then,
\[\phi(p\cdot x)=\phi(\sigma^{-1}(p)x)=\sigma^{-1}(p)\phi(x)=p\cdot \phi(x)\]
which means $^{\sigma}\phi: {^{\sigma}\mathcal{X}}\rightarrow
{^{\sigma}\mathcal{Y}}$ is a homomorphism of modules under the new
module action.

 Let
$X=(X_{i},i\in Q_{0};\;X_{\rho}: X_{s(\rho)}\rightarrow
X_{t(\rho)},\rho \in Q_{1})$ be a representation of $Q$ and
$\mathcal{X}$ with the corresponding $kQ$-module via the functor
described in Theorem~\ref{thm 2.1}, so
$\mathcal{X}=\bigoplus_{i\in Q_{0}}X_{i}$. We describe the
representation $^{\sigma}X=({^{\sigma}X}_{i},i\in
Q_{0};\;{^{\sigma}X}_{\rho}:{^{\sigma}X}_{s(\rho)}\rightarrow
{^{\sigma}X}_{t(\rho)},\rho \in Q_{1})$ corresponding to the
module $^{\sigma}\mathcal{X}$ in terms of the original
representation $X$.

\begin{proposition}\label{pro 2.2}
For a representation $X=(X_{i},i\in Q_{0};\;X_{\rho}:
X_{s(\rho)}\rightarrow X_{t(\rho)},\rho \in Q_{1})$ of a quiver
$Q$ and $\mathcal{X}$ with the corresponding $kQ$-module,  the
corresponding representation $^{\sigma}X$ of the module
$^{\sigma}\mathcal X$,
 \[^{\sigma}X=({^{\sigma}X}_{i},i\in
Q_{0},{^{\sigma}X}_{\rho}:{^{\sigma}X}_{s(\rho)}\rightarrow
{^{\sigma}X}_{t(\rho)},\rho \in Q_{1}),\] is given with
${^{\sigma}X}_{i}=X_{\sigma^{-1}(i)}$ as vector spaces, and the
map ${^{\sigma}X_{\rho}}: {^{\sigma}X}_{s(\rho)}\rightarrow
{^{\sigma}X}_{t(\rho)}$ is the same as $X_{\sigma^{-1}(\rho)}:
X_{\sigma^{-1}(s(\rho))}\rightarrow X_{\sigma^{-1}(t(\rho))}$.
\end{proposition}

\emph{Proof.}  For all $j\in Q_{0}$ and
 $\rho\in Q_{1},  x\in
{^{\sigma}X}_{s(\rho)}=X_{\sigma^{-1}(s(\rho))}$, we have that

${^{\sigma}X}_{j}=e_{j}\cdot{^{\sigma}\mathcal{X}}=\sigma^{-1}(e_{j}){^{\sigma}\mathcal{X}}=
e_{\sigma^{-1}(j)}(\bigoplus_{i\in
Q_{0}}X_{i})=X_{\sigma^{-1}(j)}$;

${^{\sigma}X}_{\rho}(x)=\rho\cdot
x=\rho\cdot(e_{s(\rho)}\cdot(\delta_{\sigma^{-1}(s(\rho))i}x)_{i\in
Q_{0}})=\sigma^{-1}(\rho)(e_{\sigma^{-1}(s(\rho))}(\delta_{\sigma^{-1}(s(\rho))i}x)_{i\in
Q_{0}})$

 $=\sigma^{-1}(\rho)(x)=X_{\sigma^{-1}(\rho)}(x)\in
X_{\sigma^{-1}(t(\rho))}={^{\sigma}X}_{t(\rho)}$. $\blacksquare$

Thus, let $\varphi=(\varphi_{i})_{i\in Q_{0}}: X\rightarrow Y$ be
the morphism between two representations, then
$^{\sigma}\varphi=({^{\sigma}\varphi_{i}})_{i\in
Q_{0}}:{^{\sigma}X}\rightarrow {^{\sigma}Y}$ satisfies
${^{\sigma}\varphi_{i}}=\varphi_{\sigma^{-1}(i)}$ as a linear map.

In this way, we obtain an additive equivalence functor $F(\sigma)$,
 with inverse $F(\sigma^{-1})$, on mod$kQ$(or say, on $Rep(Q)$), which send
$\mathcal{X}$ (or say, on $X$) to ${^{\sigma}\mathcal{X}}$ (or
say, to $^{\sigma}X$) and send $\phi:\mathcal{X}\rightarrow
\mathcal{Y}$ (or say,  $\varphi: X\rightarrow Y$) to
$^{\sigma}\phi: {^{\sigma}\mathcal{X}}\rightarrow
{^{\sigma}\mathcal{Y}}$ (or say, ${^{\sigma}\varphi}:
{^{\sigma}X}\rightarrow {^{\sigma}Y}$), and satisfies that
$F(\sigma^{r})=F(\sigma)^{r}$ for any integer $r$.

Of course,  $\mathcal{X}$ (or say, $X$) is indecomposable if and
only if ${^{\sigma}\mathcal{X}}$ (or say, ${^{\sigma}X}$)is so.

In the sequel, we always identify objects and morphisms of mod$kQ$
with the corresponding ones of Rep$Q$.

We call  a representation $X$ of a quiver $Q$ {\em isomorphically
invariant} by $\langle\sigma\rangle$ (or say,  {\em
$\langle\sigma\rangle$-equivalent }) if there is a representation
isomorphism $^{\sigma}X\cong X.$
 Equivalently, we can define an {\em
isomorphically invariant module} (or say, {\em
$\langle\sigma\rangle$-equivalent module}).

 Let $\mid
Q_0\mid=s$, $Z$ the set of all integers. Then, we can define
$\sigma: Z^s\rightarrow Z^s$ by
$\sigma(\underline{dim}X)=\underline{dim}{^{\sigma}X}$. So, a
$\langle\sigma\rangle$-equivalent representation has a dimension
vector fixed by this $\sigma$.

\begin{lemma}$(\cite{H,H2})$ \label{lem 2.3}Any indecomposable $\langle\sigma\rangle$-equivalent representation $X$ in Rep$Q$ is precisely the
representation of the form
\[X\cong Y\oplus{^{\sigma}Y}\oplus\cdots\oplus{^{\sigma^{m-1}}Y}\]
where $Y$ is an indecomposable $Q$-representation and $m\geq 1$ is
the minimal integer such that $^{\sigma^{m}}Y\cong Y.$ Moreover,
the Krull-Remak-Schmidt Theorem holds for
$\langle\sigma\rangle$-equivalent representations.
\end{lemma}

We call $X$ an {\em indecomposable
$\langle\sigma\rangle$-equivalent representation} if it is not
isomorphic to the proper direct sum of two
$\langle\sigma\rangle$-equivalent representations. Lemma~\ref{lem
2.3} means that any $\langle\sigma\rangle$-equivalent
representation is a direct sum of indecomposable
$\langle\sigma\rangle$-equivalent representations. The following
lemma tells us the relation between $m$ and $n$.

\begin{lemma} \label{lem 2.4} For any indecomposable $\langle\sigma\rangle$-equivalent $kQ$-module $Y\oplus {^{\sigma}Y}\oplus\cdots\oplus
{^{\sigma^{m-1}}Y}$, with $Y$ an indecomposable $kQ$-module and
$m$ the minimal integer such that ${^{\sigma^{m}}Y\cong Y}$. Then
$r=n/m$ is an integer.
\end{lemma}

\emph{Proof.} From ${^{\sigma^{m}}Y\cong Y}$, we have that
${^{\sigma^{km+l}}X\cong {^{\sigma^{l}}}X} $, for any non-negative
integer $k$ and  $ l\in \{0, 1,\cdots, m-1\}$. And from
${^{\sigma^{n}}X=X}$ and the minimality of $m$ such that
${^{\sigma^{m}}X\cong X}$, we have that there exists an integer $r$
such that $rm=n$. $\blacksquare$

\begin{example} Let $Q$ be the quiver $-1\stackrel{\alpha}\rightarrow 0
\stackrel{\beta}\leftarrow 1$ and $\sigma\in Aut(Q)$, which is
defined as $\sigma(e_{-1})=e_{1},
\sigma(e_{0})=e_{0},\sigma(e_{1})=e_{-1},\sigma(\alpha)=\beta,\sigma(\beta)=\alpha.$
All indecomposable $kQ$-representations are:

\[L_{-1}: k\rightarrow 0 \leftarrow 0;\ \ \ \ \ \
L_{0}: 0\rightarrow k \leftarrow
 0;\ \ \ \ \ \ L_{1}: 0\rightarrow 0 \leftarrow k;\]
\[L_{-1 0}: k\stackrel{1}\rightarrow k \leftarrow 0;\ \ \ \ \ \ L_{0
1}: 0\rightarrow k \stackrel{1}\leftarrow k;\ \ \ \ \ \ L_{1 0 1}:
k\stackrel{1}\rightarrow k \stackrel{1}\leftarrow
 k.\]

By Proposition~\ref{pro 2.2},  we have
$^{\sigma}L_{-1}=L_{1},{^{\sigma}L_{1}}=L_{-1},{^{\sigma}L_{-10}}=L_{01},
 {^{\sigma}L_{01}}=L_{-10},{^{\sigma}L_{0}}=L_{0},{^{\sigma}L_{101}}=L_{101}.$
By definitions, all indecomposable
$\langle\sigma\rangle$-equivalent representations are:
\[L_{1}\oplus L_{-1};\ \ \ \ L_{01}\oplus L_{-10};\ \ \ \
L_{0};\ \ \ \ L_{101}.\]

\end{example}

\section{Structure of modules over a skew group algebra from a cyclic group }

  In this section, we denote by $(Q,\sigma)$ a  fixed  connected finite quiver $Q$ without oriented cycle and a quiver automorphism $\sigma\in AutQ$
 of order $n$.  Then, we have a cyclic group $\langle\sigma\rangle$ of order $n$ and a skew group algebra $kQ\#k\langle\sigma\rangle$. In the sequel, we will always assume that $k$ is an algebraically
  closed  field  with the characteristic of not dividing $n$.

  The following lemma can be found in \cite{H}, but we still give its
  proof because it is useful for our discussion.

  \begin{lemma}$(\cite{H,H2})$ \label{lem 3.1}
Every module $X$ of the skew group
  algebra $kQ\#k\langle\sigma\rangle$ is a $\langle\sigma\rangle$-equivalent $kQ$-module.
  \end{lemma}
  \emph{Proof.} In order to show $X$ to be a $\langle\sigma\rangle$-equivalent $kQ$-module, we
  only need to prove $X\cong {^{\sigma}X}$ as $kQ$-modules. Define $f:{^{\sigma}X}\rightarrow
  X$ such that $f(x)=\sigma x$, for all $x\in X$. It is
  well-defined since $X$ is a $kQ\#k\langle\sigma\rangle$-module.
 We have that for any $ p\in kQ$,
  \[f(p\cdot x)=\sigma(p\cdot x)=\sigma(\sigma^{-1}(p)x)=\sigma(\sigma^{-1}(p))x=(p\#\sigma)x=p(\sigma
  x)=pf(x),\]
  which means that $f$ is a
  $kQ$-module homomorphism. Moreover, $f$ is an isomorphism with inverse $f^{-1}: X\rightarrow
  {^{\sigma}X}$ such that $f^{-1}(x)=\sigma^{-1}x$ for all $x\in
  X$. $\blacksquare$

By Lemma~\ref{lem 3.1} and Lemma~\ref{lem 2.3}, we have for any
$kQ\#k\langle\sigma\rangle$-module $X$, $X\cong_{i=1}^{s} Y_{i}$,
where $Y_{i}\cong
X_i\oplus{^{\sigma}X_i}\oplus\cdots\oplus{^{\sigma^{m_i-1}}X_i}$
with $X_{i}$ an indecomposable $kQ$-module and $m_{i}$ a minimal
positive integer such that ${^{\sigma^{m_{i}}}X_{i}}\cong
  X_{i}$. Hence we have the following $kQ$-isomorphism, say $g$:
\[X\stackrel{g}\cong \bigoplus_{i=1}^{s}\bigoplus_{j=0}^{m_{i}-1}
   {^{\sigma^{j}}X_{i}}.\]
Then we can define the $kQ\#k\langle\sigma\rangle$-module
structure on $\bigoplus_{i=1}^{s}\bigoplus_{j=0}^{m_{i}-1}
   {^{\sigma^{j}}X_{i}}$ through $g$.
In fact, we define $\sigma'$s action on
$\bigoplus_{i=1}^{s}\bigoplus_{j=0}^{m_{i}-1}
   {^{\sigma^{j}}X_{i}}$ by  \[\sigma
y=g(\sigma g^{-1}(y)),\] and $(p\sigma^l)y=p(\sigma^ly)$  for any
$p\in kQ$, $y\in \bigoplus_{i=1}^{s}\bigoplus_{j=0}^{m_{i}-1}
   {^{\sigma^{j}}X_{i}}, l=0,\cdots, n-1$.

  The action makes $\bigoplus_{i=1}^{s}\bigoplus_{j=0}^{m_{i}-1}
   {^{\sigma^{j}}X_{i}}$ a $kQ\#k\langle\sigma\rangle$-module since by definition of the smash product, $\sigma p=\sigma(p)\#\sigma$,
then
\begin{eqnarray*}\sigma(py)&=&g(\sigma(g^{-1}(py)))=g(\sigma(pg^{-1}(y)))=g(\sigma(p)(\sigma(g^{-1}y)))\\
   &=&\sigma(p)(g(\sigma(g^{-1}y)))=\sigma(p)(\sigma y)=(\sigma(p)\sigma)y\\
   &=&(\sigma
   p)y.
   \end{eqnarray*}
     Moreover, $g$ is a $kQ\#k\langle\sigma\rangle$-module
   homomorphism
via $g(\sigma x)=g(\sigma(g^{-1}(g(x))))=\sigma(g(x))$  for any
$x\in X$.

Considering the restriction of  $\sigma$ on each indecomposable
$kQ$-module $X_i$, we have:
  \begin{corollary} \label{cor 3.2} With the above notations, ${^{\sigma^{j}}X_{i}}\cong
  \sigma^{j}X_{i}$ as $kQ$-modules for any $j\in\{1,\cdots,m_{i}\},\;i\in \{1,\cdots,s\}$.
  \end{corollary}
  \emph{Proof.} Define $f: {^{\sigma^{j}}X_{i}}\rightarrow
  \sigma^{j}X_{i}$
  by $f(x)=\sigma^{j}x=g(\sigma^{j}(g^{-1}(x)))$ for any $x\in
  {^{\sigma^{j}}X_{i}}$.
Then, for any $p\in kQ$,
\begin{eqnarray*}f(p\cdot x)&=&\sigma^{j}(p\cdot x)=
  g(\sigma^{j}(g^{-1}(p\cdot
  x)))=g(\sigma^{j}(g^{-1}(\sigma^{-j}(p)x)))\\
  &=&g(g^{-1}(\sigma^{j}(\sigma^{-j}(p)x)))
  =\sigma^{j}(\sigma^{-j}(p)x)=(p\#\sigma^{j})x=p(\sigma^{j}(x))\\&=&p(f(x)),\end{eqnarray*}
 thus  $f$ is a $kQ$-module isomorphism with inverse $f^{-1}$ satisfying
 $f^{-1}(y)=\sigma^{-j}y$ for any $y\in\sigma^jX_i$.  $\blacksquare$

By corollary~\ref{cor 3.2}, for any
$kQ\#k\langle\sigma\rangle$-module $X$, $X\cong
\oplus_{i=1}^{s}Y_{i}$ with $Y_{i}$ as an indecomposable
$\langle\sigma\rangle$-equivalent $kQ$-module, $\sigma'$s action
may be closed in each $Y_{i}$. Interestingly, whether $\sigma'$s
action is closed in each $Y_{i}$ is the same thing as whether any
indecomposable $kQ\#k\langle\sigma\rangle$-module is an
indecomposable $\langle\sigma\rangle$-equivalent $kQ$-module.

   From any indecomposable $kQ$-module $X$, we
 get an indecomposable $\langle\sigma\rangle$-equivalent $kQ$-module
  $\oplus_{i=1}^{m-1}{^{\sigma^{i}}X}$, and then we have several questions to
  solve:

\emph{\bf Question 1.}  Is
  it possible to endow $\oplus_{i=1}^{m-1}{^{\sigma^{i}}X}$ with
  an induced $kQ\#k\langle\sigma\rangle$-module structure?

 \emph{\bf Question 2.} If
 the $kQ$-module $\bigoplus_{i=1}^{m-1}{^{\sigma^{i}}X}$ can be endowed with some
  $kQ\#k\langle\sigma\rangle$-module structures, how many non-isomorphic classes of such induced
   $kQ\#k\langle\sigma\rangle$-modules exist?

\section{Construction of indecomposable
$kQ\#k\langle\sigma\rangle$-modules
 from indecomposable $kQ$-modules
  }

It is known in \cite{H,H2} that any
$\langle\sigma\rangle$-equivalent $kQ$-module $X$ induces a
$kQ\#k\langle\sigma\rangle$-module. For   completeness, we give
its proof below:

\begin{proposition} $(\cite{H,H2})$ \label{pro 4.1} Let $X$ be a $\langle\sigma\rangle$-equivalent $kQ$-module. Then
there exists an isomorphism $\phi: {^{\sigma}X}\rightarrow X$ such
that $\phi^{n}\;=\;\phi {\;^{\sigma}\phi}\cdots
{\;^{\sigma^{n-1}}\phi}$ of $X$ is the identity.
\end{proposition}

\begin{theorem}$(\cite{H,H2})$  \label{thm 4.2} Let $X$ be a $\langle\sigma\rangle$-equivalent $kQ$-module. Then
$X$  has an induced $kQ\#k\langle\sigma\rangle$-module structure if
we define $\sigma^{i}(x)=\phi^{i}(x)$.
\end{theorem}
\emph{Proof.} Since for $p\#\sigma^{i},q\#\sigma^{j}\in
kQ\#k\langle\sigma\rangle, x\in X$,
\begin{eqnarray*}(p\#\sigma^{i})((q\#\sigma^{j})(x))&=&(p\#\sigma^{i})(q\phi^{j}(x))=p\phi^{i}(q\phi^{j}(x))=p\phi^{i-1}\phi(\sigma(q)\cdot
\phi^{j}(x))\\&=&p\phi^{i-1}(\sigma(q)\phi^{j+1}(x))=\cdots=p\sigma^{i}(q)(\phi^{i+j}(x))\\&=&(p\sigma^{i}(q)\#\sigma^{i+j})(x)=(p\#\sigma^{i})(q\#\sigma^{j})(x).
\;\;\;\;\;\;\blacksquare
\end{eqnarray*}

This result means Question 1 can be answered affirmatively.

The induced $kQ\#k\langle\sigma\rangle$-module constructed in
Theorem \ref{thm 4.2} will be used accordingly in our conclusions
below. So, we will call such induced module a {\em canonical
induced $kQ\#k\langle\sigma\rangle$-module}.

Again, we recall two lemmas, which will be used in the next two
sections.

\begin{lemma}\label{le3.7}$(\cite{RR})$ Let $X,Y$ be indecomposable $kQ$-modules, and $G$ be a subgroup of the $k$-automorphism group of
 $kQ$. Then:
 \begin{itemize}

  \item[$(i)$]  $(kQ\#kG)\otimes_{kQ}X\cong \bigoplus_{g\in G}{^{g}X}$ as
  $kQ$-modules;

   \item[$(ii)$]  $(kQ\#kG)\otimes_{kQ}X\cong
  (kQ\#kG)\otimes_{kQ}Y$ if and only if $Y\cong {^{g}X}$ for
  some $g\in G$;

   \item[$(iii)$]  The number of summands in the decomposition of
  $(kQ\#kG)\otimes_{kQ}X$ into a direct sum of indecomposables is at most the order of $H$, where $H=\{g\in G, {^{g}X}\cong
  X\}$;

   \item[$(iv)$]  If $G$ is cyclic of order $n$ and $X\cong {^{g}X}$
  for all $g\in G$, then $(kQ\#kG)\otimes_{kQ}X$ has exactly $n$
  summands;

   \item[$(v)$] If $H=\{g\in G, {^{g}X}\cong X\}$ is cyclic of order
  $m$, then $(kQ\#kG)\otimes_{kQ}X$ has exactly $m$ summands.

  \end{itemize}

  \end{lemma}

 \begin{lemma}\label{le3.8} $(\cite{L,LZ})$  Let $H$ be a finite dimensional semisimple Hopf
  algebra and $A$ be a finite dimensional $H$-module algebra. Then, for any $A\#H$-module $X$, it holds that
$X\mid(A\#H)\otimes_{A}X$, that is, $X$ is a direct summand of
$(A\#H)\otimes_{A}X$ as an $A\#H$-module.
  \end{lemma}

\subsection{Induction of indecomposable $kQ\#k\langle\sigma^{m}\rangle$-modules from an indecomposable $kQ$-module $X$ with minimal $m$ satisfying
 $^{\sigma^{m}}X\cong X$}

 For any indecomposable $kQ$-module $X$ with minimal $m$ such that
${^{\sigma^{m}}X}\cong X$, let $L=\{g\in \langle\sigma\rangle \mid
{^{g}X}\cong X\}$, then $L=\langle\sigma^{m}\rangle$, a cyclic
group generated by $\sigma^{m}$ by Lemma~\ref{lem 2.4}. Since
$kL\cong k\langle\sigma^{m}\rangle$ is a semisimple group algebra,
we have
\begin{equation} \label{eq 1}
k\langle\sigma^{m}\rangle\cong\bigoplus_{i=1}^{r=n/m} L_{i}
\end{equation}
 as a
$k\langle\sigma^{m}\rangle$-module, where $L_{i}$ is isomorphic to
$k$ as a vector space, and ${\sigma^{m}}'$s action is
$\sigma^{m}(1)=\zeta^{i}$, where $\zeta$ is the $r$-th primitive
root of 1. Moreover, $L_{i}\ncong L_{j}$ as
$k\langle\sigma^{m}\rangle$-modules, if $i\neq j$.

Before answering Question 2, we introduce the following question:

\emph{\bf Question 3.} For any indecomposable $kQ$-module $X$ with
minimal $m$ such that ${^{\sigma^{m}}X\cong X}$, how many
non-isomorphic indecomposable
$kQ\#k\langle\sigma^{m}\rangle$-modules can be induced from $X$?

Since $X$ is a $\langle\sigma^{m}\rangle$-equivalent $kQ$-module,
by Proposition~\ref{pro 4.1} and Theorem~\ref{thm 4.2}, there
exists $\phi_{X}$, such that  $X$ induces a
$kQ\#k\langle\sigma^{m}\rangle$-module structure. Using the
$kQ\#k\langle\sigma^{m}\rangle$-module structure on $X$, we can
define $kQ\#k\langle\sigma^{m}\rangle$-module structure on
$L_{i}\otimes_{k}X$, for any $i\in\{1,\cdots,r\}$. In fact, for
any $ i,j\in \{1,2,\cdots,r\}, 1\otimes x\in L_{i}\otimes_{k}
X,p\in kQ,$ define
\[p\#\sigma^{mj}(1\otimes x)=\sigma^{mj}(1)\otimes p\#l(x)=\zeta^{ij}\otimes p\#l(x), \]
 where the action $p\#l(x)$ is inherited from the canonical
 induced
 $kQ\#k\langle\sigma^{m}\rangle$-module structure. Then

\begin{lemma} \label{lem 4.3} With the above notations, we have:
\begin{itemize}

\item[$(i)$] $L_{i}\otimes_{k}X\cong X$ as $kQ$-modules for any
$i\in \{1,2,\cdots,r\}$;

\item[$(ii)$] $L_{i}\otimes_{k}X$ is an indecomposable
$kQ\#k\langle\sigma^{m}\rangle$-module for any $i\in
\{1,2,\cdots,r\}$;

\item[$(iii)$] $L_{i}\otimes_{k}X\ncong L_{j}\otimes_{k}X$ as
$kQ\#k\langle\sigma^{m}\rangle$-modules, if $i\neq j$.

\end{itemize}

\end{lemma}
\emph{Proof.} $(i)$ Define $f: X\rightarrow L_{i}\otimes_{k} X$ by
$x\longmapsto 1\otimes x$, then $f$ is a $kQ$-module homomorphism
since $f(p(x))=1\otimes p(x)=p(1\otimes x)=pf(x)$, for any $ p\in
kQ, x\in X$. Obviously, $f$ is bijective.

$(ii)$ $\;L_{i}\otimes_{k}X$ is
$kQ\#k\langle\sigma^{m}\rangle$-module, since for  any $
p_{1},p_{2}\in kQ, j_{1},j_{2}\in \{1,2,\cdots, r\}, x\in X$,
\begin{eqnarray*}
(1\#1)(1\otimes x)&=&\sigma^{mr}(1)\otimes 1\#1(x)=\zeta^{ir}\otimes x=1\otimes x,\\
((p_{2}\#\sigma^{mj_{2}})(p_{1}\#\sigma^{mj_{1}}))(1\otimes
x)&=&(p_{2}\sigma^{mj_{2}}(p_{1})\#\sigma^{mj_{2}}\sigma^{mj_{1}})(1\otimes
x)\\&=&(p_{2}\sigma^{mj_{2}}(p_{1})\#\sigma^{m(j_{2}+j_{1})})(1\otimes
x)\\&=&\sigma^{m(j_{2}+j_{1})}(1)\otimes
(p_{2}\sigma^{mj_{2}}(p_{1})\#\sigma^{mj_{2}}\sigma^{mj_{1}})(x)\\
&=&\zeta^{i(j_{2}+j_{1})}(1)\otimes
(p_{2}\sigma^{mj_{2}}(p_{1})\#\sigma^{mj_{2}}\sigma^{mj_{1}})(x)\\&=&
\zeta^{ij_{2}}\zeta^{ij_{1}}(1)\otimes
(p_{2}\#\sigma^{mj_{2}})((p_{1}\#\sigma^{mj_{1}})(x))\\&=&p_{2}\#\sigma^{mj_{2}}(\zeta^{ij_{1}}(1)\otimes p_{1}\#\sigma^{mj_{1}}(x))\\
&=& p_{2}\#\sigma^{mj_{2}}(p_{1}\#\sigma^{mj_{1}}(1\otimes x)).
\end{eqnarray*}

And $L_{i}\otimes_{k}X$ is an indecomposable
$kQ\#k\langle\sigma^{m}\rangle$-module since  it is an
indecomposable $kQ$-module by $(i)$.

Before giving a proof of $(iii)$, we need to perform some
preparations as follows.

\begin{lemma} \label{lem 4.4} For any $i\in \{1,2,\cdots,r\},$  $Hom_{kQ}(X, L_{i}\otimes_{k} X)\cong
L_{i}\otimes_{k} End_{kQ}(X)$ as
$kQ\#k\langle\sigma^{m}\rangle$-modules.
\end{lemma}

\emph{Proof.} The $kQ\#k\langle\sigma^{m}\rangle$-module structure
of $Hom_{kQ}(X, L_{i}\otimes_{k} X)$ is given by
$$(p\#l(f))(x)=(p\#l)f(x), \forall  f\in Hom_{kQ}(X,
L_{i}\otimes_{k} X), p\#l\in kQ\#k\langle\sigma^{m}\rangle, x\in
X;$$

The $kQ\#k\langle\sigma^{m}\rangle$-module structure of
$End_{kQ}(X)$ is given by
$$(p\#l(f))(x)=(p\#l)f(x),\forall f\in End_{kQ}(X), p\#l\in
kQ\#k\langle\sigma^{m}\rangle, x\in X;$$

  The $kQ\#k\langle\sigma^{m}\rangle$-module
structure of $L_{i}\otimes_{k} End_{kQ}(X)$ is given by
$$(p\#l)(1\#f)=l(1)\#(p\#l)(f),\forall f\in End_{kQ}(X), p\#l\in
kQ\#k\langle\sigma^{m}\rangle.$$

Define $F: Hom_{kQ}(X, L_{i}\otimes_{k} X)\rightarrow L_{i}\otimes
End_{kQ}(X)$ by $$F(f)=1\otimes \overline{f}, \forall f\in
Hom_{kQ}(X, L_{i}\otimes_{k} X),$$ where $\overline{f}$ is defined
 by $\overline{f}(x)=k_{f}x_{f}$, if $f(x)=k_{f}\otimes
x_{f},\forall x\in X$. Since for any $p\in kQ, x\in X$,
$f(px)=pf(x)=p(k_{f}\otimes x_{f})=k_{f}\otimes p(x_{f})=1\otimes
p(k_{f}x_{f})$, then
$\overline{f}(px)=k_{f}p(x_{f})=p(k_{f}x_{f})=p\overline{f}(x)$,
i.e., $\overline{f}\in End_{kQ}(X)$, which means $F$ is
well-defined.

 Show $F$ is a $kQ\#k\langle\sigma\rangle$-module
homomorphism.  Since for any $p\in kQ, f\in Hom_{kQ}(X,
L_{i}\otimes_{k} X) $,  $(pf)(x)=p(f(x))=p(k_{f}\otimes
x_{f})=k_{f}\otimes p(x)$, then
   $F(pf)=1\otimes \overline{pf}=1\otimes
p\overline{f}=pF(f) $, which means $F$ is a $kQ$-module
homomorphism. And since
$(\sigma^{m}f)(x)=\sigma^{m}(f(x))=\sigma^{m}(k_{f}\otimes
x_{f})=\zeta^{i}k_{f}\otimes \sigma^{m} x_{f}$, then
$F(\sigma^{m}f)=1\otimes \overline{\sigma^{m}f}=1\otimes \zeta^{i}
\sigma^{m}\overline{f}=\zeta^{i}\otimes
\sigma^{m}\overline{f}=\sigma^{m}(1\otimes \overline{f})$, which
means $F$ is a $k\langle\sigma^{m}\rangle$-module homomorphism.

Finally, $F$ is a $kQ\#k\langle\sigma^{m}\rangle$-module
isomorphism, since $F$ is injective and  $$dim_k Hom_{kQ}(X,
L_{i}\otimes_{k} X)=dim_k L_{i}\otimes_{k} End_{kQ}(X).
\;\;\;\blacksquare$$

Now we go back to the proof of $(iii)$:

 Otherwise,
$L_{i}\otimes_{k}X\cong L_{j}\otimes_{k} X$ as
$kQ\#k\langle\sigma^{m}\rangle$-modules, for some $i\neq j$,  then
by Lemma~\ref{lem 4.4}, $L_{i}\otimes_{k}End_{kQ}(X)\cong
L_{j}\otimes_{k}End_{kQ}(X)$ as
$kQ\#k\langle\sigma^{m}\rangle$-modules.  Since $End_{kQ}(X)$ is a
local ring, $k$ is algebraically closed, we have
$End_{kQ}(X)/radEnd_{kQ}(X)\cong k$ as algebras. Since
$radEnd_{kQ}(X)$ is closed under
$k\langle\sigma^{m}\rangle$-module structure, we have
$L_{i}\otimes_{k}End_{kQ}(X)/radEnd_{kQ}(X)\cong
L_{j}\otimes_{k}End_{kQ}(X)/radEnd_{kQ}(X)$, which induces
$L_{i}\cong L_{j}$ as $k\langle\sigma^{m}\rangle$-modules, with
contradiction to $i\neq j$.  $\blacksquare$

\begin{theorem}\label{Thm 4.7} Let $X$ be an indecomposable $kQ$-module with $m$ minimal such that $\sigma^{m}X\cong X$,  $L_{i}$ is defined in the isomorphism relation~$(\ref{eq 1})$, for any $i\in \{1,2,\cdots, r\}$.
 Then the following statements hold:
\begin{itemize}
\item[$(i)$] $(kQ\#k\langle\sigma^{m}\rangle)\otimes_{kQ} X$ is
isomorphic to the direct sum of $r$ non-isomorphic
$kQ\#k\langle\sigma^{m}\rangle$-modules, that is,
$(kQ\#k\langle\sigma^{m}\rangle)\otimes_{kQ} X\cong
\bigoplus_{i=1}^{r}L_{i}\otimes_{k}X$ as
$kQ\#k\langle\sigma^{m}\rangle$-modules;

\item[$(ii)$] For any $kQ\#k\langle\sigma^{m}\rangle$-module $Y$,
if $Y\cong X$ as $kQ$-modules, then there exists a unique $i\in
\{1,2,\cdots, r\}$, such that $Y\cong L_{i}\otimes_{k}X$. That is,
there are $r$ non-isomorphic
$kQ\#k\langle\sigma^{m}\rangle$-modules induced from $X$.
\end{itemize}
\end{theorem}
{\em Proof.} $(i)$ For any $i\in \{1,2,\cdots,r\}$, since
$L_{i}\otimes_{k}X\mid kQ\#k\langle\sigma^{m}\rangle\otimes_{kQ}
  (L_{i}\otimes_{k}X)$ by Lemma~\ref{le3.8},  $kQ\#k\langle\sigma^{m}\rangle\otimes_{kQ}
  (L_{i}\otimes_{k}X)\cong
  kQ\#k\langle\sigma^{m}\rangle\otimes_{kQ}X$ by Lemma~\ref{lem 4.3}$(i)$, then we have $L_{i}\otimes_{k}X\mid
kQ\#k\langle\sigma^{m}\rangle\otimes_{kQ}
  X$. Then by by Lemma~\ref{lem
4.3}$(iii)$, if $i\neq j,$ $L_{i}\otimes_{k}X\ncong
L_{j}\otimes_{k}X,$ then $(\oplus_{i=1}^{r}L_{i}\otimes_{k}X)\mid
kQ\#k\langle\sigma^{m}\rangle\otimes_{kQ}
  X$ by Krull-Schmidt Theorem. And $kQ\#k\langle\sigma^{m}\rangle\otimes_{kQ} X\cong
\bigoplus_{i=1}^{r}L_{i}\otimes_{k}X$ since
$kQ\#k\langle\sigma^{m}\rangle\otimes_{kQ}X$ has exactly $r$
  indecomposable summands by Lemma~\ref{le3.7}$(v)$.

$(ii)$ For a $kQ\#k\langle\sigma^{m}\rangle$-module $Y$,  $Y\cong
X$ as $kQ$-modules, then $Y$ is an indecomposable
$kQ\#k\langle\sigma^{m}\rangle$-module and by Lemma~\ref{le3.8},
$Y\mid kQ\#k\langle\sigma^{m}\rangle\otimes_{kQ}Y\cong kQ
\#k\langle\sigma^{m}\rangle\otimes_{kQ}X$.  Then by $(i)$,
Lemma~\ref{lem 4.3}$(iii)$ and the Krull-Schmidt Theorem, there
exists a unique $i\in\{1,2,\cdots, r\}$, such that $Y\cong
L_{i}\otimes_{k}X$. $\blacksquare$

\subsection{Induction of indecomposable
$kQ\#k\langle\sigma\rangle$-modules
 from an indecomposable $\langle\sigma\rangle$-equivalent
 $kQ$-module
}

In this section, we are ready to answer Question 2.

\begin{theorem}\label{thm 4.8} Let $X$ be an indecomposable $kQ$-module with $m$ minimal such that $\sigma^{m}X\cong X$, $L_{i}$
as defined in the isomorphism relation~$(\ref{eq 1})$, for any
$i\in\{1,2,\cdots, r\}$.
 Then the following statements hold:

\begin{itemize}
\item[$(i)$]
$(kQ\#k\langle\sigma\rangle)\otimes_{kQ\#k\langle\sigma^{m}\rangle}(L_{i}\otimes_{k}X)\cong
X\oplus {^{\sigma}X}\oplus\cdots \oplus{^{\sigma^{m-1}}X}$ as
$kQ$-modules;

\item[$(ii)$]
$(kQ\#k\langle\sigma\rangle)\otimes_{kQ\#k\langle\sigma^{m}\rangle}(L_{i}\otimes_{k}X)$
is an indecomposable $kQ\#k\langle\sigma\rangle$-module;

\item[$(iii)$]
$(kQ\#k\langle\sigma\rangle)\otimes_{kQ\#k\langle\sigma^{m}\rangle}(L_{i}\otimes_{k}X)\ncong
(kQ\#k\langle\sigma\rangle)\otimes_{kQ\#k\langle\sigma^{m}\rangle}(L_{j}\otimes_{k}X)$
as $kQ\#k\langle\sigma\rangle$-module, if $i\neq j$;

\item[$(iv)$] $(kQ\#k\langle\sigma\rangle)\otimes_{kQ} X$ is
isomorphic to the direct sum of $r$ non-isomorphic
$kQ\#k\langle\sigma\rangle$-modules, that is,
$(kQ\#k\langle\sigma\rangle)\otimes_{kQ}X\cong
\oplus_{i=1}^{r}(kQ\#k\langle\sigma\rangle)\otimes_{kQ\#k\langle\sigma^{m}\rangle}(L_{i}\otimes_{k}X)$
as $kQ\#k\langle\sigma\rangle$-modules;

\item[$(v)$] For any $kQ\#k\langle\sigma\rangle$-module $Y$, if
$Y\cong X\oplus {^{\sigma}X}\oplus\cdots \oplus{^{\sigma^{m-1}}X}$
as $kQ$-modules, then there exists a unique $i\in \{1,2,\cdots,
r\}$, such that $Y\cong
(kQ\#k\langle\sigma\rangle)\otimes_{kQ\#k\langle\sigma^{m}\rangle}(L_{i}\otimes_{k}X)$.
That is, there are $r$ non-isomorphic
$kQ\#k\langle\sigma\rangle$-modules induced from the
indecomposable $\langle\sigma\rangle$-equivalent $kQ$-modules
$X\oplus {^{\sigma}X}\oplus\cdots \oplus{^{\sigma^{m-1}}X}$.
\end{itemize}
\end{theorem}
\emph{Proof.} $(i)$
$(kQ\#k\langle\sigma\rangle)\otimes_{kQ\#k\langle\sigma^{m}\rangle}(L_{i}\otimes_{k}X)=1\otimes
X\oplus \sigma\otimes X\cdots \oplus \sigma^{m-1}\otimes X\cong
X\oplus {^{\sigma}X}\oplus\cdots \oplus{^{\sigma^{m-1}}X}$ as
$kQ$-modules since $ {^{\sigma^{j}}X}\cong \sigma^{j}\otimes X$ as
$kQ$-modules. In fact, define $f: {^{\sigma^{j}}X}\rightarrow
\sigma^{j}\otimes X$ by $f(x)=\sigma^{j}\otimes x,$ for any $x\in
X$. Then $f$ is bijection, and $f$ is a $kQ$-module homomorphism
since $f(p\cdot x)=\sigma^{j}\otimes
\sigma^{-j}(p)(x)=\sigma^{j}(\sigma^{-j}(p))\otimes
x=p(\sigma^{j}\otimes x)=pf(x), \forall p\in kQ, x\in X.$

$(ii)$
$(kQ\#k\langle\sigma\rangle)\otimes_{kQ\#k\langle\sigma^{m}\rangle}(L_{i}\otimes_{k}X)$
is an indecomposable $kQ\#k\langle\sigma\rangle$-module, since
$kQ\#k\langle\sigma\rangle\otimes_{kQ\#k\langle\sigma^{m}\rangle}(L_{i}\otimes_{k}X)$
is an indecomposable $\langle\sigma\rangle$-equivalent $kQ$-module
by $(i)$ and Lemma~\ref{lem 3.1}.

$(iii)$ For any $i\in\{1,2,\cdots,r\}$,
$(kQ\#k\langle\sigma\rangle)\otimes_{kQ\#k\langle\sigma^{m}\rangle}(L_{i}\otimes_{k}X)=1\otimes
L_{i}\otimes_{k}X\oplus \sigma\otimes L_{i}\otimes_{k}X\cdots
\oplus \sigma^{m-1}\otimes L_{i}\otimes_{k}X$ as
$kQ\#k\langle\sigma^{m}\rangle$-modules. Otherwise, if  $i\neq j,$
$kQ\#k\langle\sigma\rangle\otimes_{kQ\#k\langle\sigma^{m}\rangle}(L_{i}\otimes_{k}X)\cong
kQ\#k\langle\sigma\rangle\otimes_{kQ\#k\langle\sigma^{m}\rangle}(L_{j}\otimes_{k}X)$
as $kQ\#k\langle\sigma\rangle$-modules. Then $1\otimes
L_{i}\otimes_{k}X\oplus \sigma\otimes L_{i}\otimes_{k}X\cdots
\oplus \sigma^{m-1}\otimes L_{i}\otimes_{k}X\cong 1\otimes
L_{j}\otimes_{k}X\oplus \sigma\otimes L_{j}\otimes_{k}X\cdots
\oplus \sigma^{m-1}\otimes L_{j}\otimes_{k}X$ as
$kQ\#k\langle\sigma^{m}\rangle$-modules, which is a contradiction
since for $1\leq s\leq m-1$,  $1\otimes L_{i}\otimes_{k}X\cong
X\ncong {^{\sigma^{s}}X}\cong \sigma^{s}\otimes L_{j}\otimes_{k}X$
as $kQ$-modules  and $1\otimes L_{i}\otimes_{k}X\cong
L_{i}\otimes_{k}X\ncong L_{j}\otimes_{k}X\cong 1\otimes
L_{j}\otimes_{k}X$ as $kQ\#k\langle\sigma^{m}\rangle$-modules.

$(iv)$ For any $i\in\{1,2,\cdots,r\}$, we have
$kQ\#k\langle\sigma\rangle\otimes_{kQ\#k\langle\sigma^{m}\rangle}(L_{i}\otimes_{k}X)\mid
kQ\#k\langle\sigma\rangle\otimes_{kQ}kQ\#k\langle\sigma\rangle\otimes_{kQ\#k\langle\sigma^{m}\rangle}(L_{i}\otimes_{k}X)$
by Lemma~\ref{le3.8},
$kQ\#k\langle\sigma\rangle\otimes_{kQ\#k\langle\sigma^{m}\rangle}(L_{i}\otimes_{k}X)\mid
kQ\#k\langle\sigma\rangle\otimes_{kQ}(X\oplus
{^{\sigma}X}\oplus\cdots \oplus{^{\sigma^{m-1}}X})$ by $(i)$, and
$kQ\#k\langle\sigma\rangle\otimes_{kQ\#k\langle\sigma^{m}\rangle}(L_{i}\otimes_{k}X)\mid
kQ\#k\langle\sigma\rangle\otimes_{kQ}X$ by $(ii)$ and
Lemma~\ref{le3.7}$(ii)$. By $(iii)$ and the Krull-Schmidt Theorem,
we have
$(\oplus_{i=1}^{r}kQ\#k\langle\sigma\rangle\otimes_{kQ\#k\langle\sigma^{m}\rangle}(L_{i}\otimes_{k}X))\mid
kQ\#k\langle\sigma\rangle\otimes_{kQ}X$. Additionally
$kQ\#k\langle\sigma\rangle\otimes_{k}X\cong
\oplus_{i=1}^{r}kQ\#k\langle\sigma\rangle\otimes_{kQ\#k\langle\sigma^{m}\rangle}(L_{i}\otimes_{k}X)$
since $kQ\#k\langle\sigma\rangle\otimes_{kQ}X$ has exactly $r$
indecomposable summands by Lemma~\ref{le3.7}$(v)$.

$(v)$ For a $kQ\#k\langle\sigma\rangle$-module $Y$,  $Y\cong
X\oplus {^{\sigma}X}\oplus\cdots \oplus{^{\sigma^{m-1}}X}$ as
$kQ$-modules, then $Y$ is an indecomposable
$kQ\#k\langle\sigma\rangle$-module and by Lemma~\ref{le3.8},
 $Y\mid
kQ\#k\langle\sigma\rangle\otimes_{kQ}Y\cong
kQ\#k\langle\sigma\rangle\otimes_{kQ}(X\oplus
{^{\sigma}X}\oplus\cdots \oplus{^{\sigma^{m-1}}X})$, then by
$(iv)$, Lemma~\ref{le3.7}$(ii)$ and the Krull-Schmidt Theorem,
there exists a unique $i\in \{1,2,\cdots,r\}$, such that $Y\cong
kQ\#k\langle\sigma\rangle\otimes_{kQ\#k\langle\sigma^{m}\rangle}L_{i}\otimes_{k}X$.
$\blacksquare$

\begin{theorem} \label{thm 4.9} Any indecomposable
$kQ\#k\langle\sigma\rangle$-module is an indecomposable
$\langle\sigma\rangle$-equivalent $kQ$-module. Conversely, for any
indecomposable $\langle\sigma\rangle$-equivalent $kQ$-module, the
corresponding canonical induced $kQ\#k\langle\sigma\rangle$-module
is indecomposable.
\end{theorem}
\emph{Proof.} Given any indecomposable
$kQ\#k\langle\sigma\rangle$-module $X$,  by Lemma~\ref{lem 2.3},
$X\cong \oplus_{j=1}^{s}X_{j}$ ,  $X_{j}\cong Y_{j}\oplus
{^{\sigma}Y_{j}}\oplus \cdots \oplus {^{\sigma^{m_{j}-1}}Y_{j}}$
with $Y_{j}$ an indecomposable $kQ$-module and $m_{j}$  minimal
such that ${^{\sigma^{m_{j}}}Y_{j}\cong Y_{j}}$.  By
Lemma~\ref{le3.8}, we have $X\mid
kQ\#k\langle\sigma\rangle\otimes_{kQ} X\cong
\oplus_{j=1}^{s}\oplus_{k=0}^{m_{j}-1}kQ\#k\langle\sigma\rangle\otimes_{kQ}
{^{\sigma^{k}}Y_{j}}$, then by Lemma~\ref{le3.7}$(ii)$ and the
Krull-Schmidt Theorem, there exists $j$, such that $X\mid
kQ\#k\langle\sigma\rangle\otimes_{kQ}Y_{j}$. Thus by
Theorem~\ref{thm 4.8}, we have $X\cong Y_{j}\oplus
{^{\sigma}Y_{j}}\oplus\cdots\oplus {^{\sigma^{m_{j}-1}}Y_{j}}$ as
$kQ$-modules, that is, $X$ is an indecomposable
$\langle\sigma\rangle$-equivalent $kQ$-modules.

Conversely, since any $kQ\#k\langle\sigma\rangle$-module is a
$\langle\sigma\rangle$-equivalent $kQ$-module. $\blacksquare$\\

 According to Theorem~\ref{thm 4.8} and Theorem~\ref{thm 4.9}, our main purpose has been
carried out, that is, all indecomposable
$kQ\#k\langle\sigma\rangle$-modules can be constructed from
indecomposable $kQ$-modules as follows:

$(I)$ For a fixed indecomposable $kQ$-module $X$, write $m$ to be
the minimal positive integer satisfying $\sigma^{m}X\cong X$.
  On the indecomposable $\langle\sigma\rangle$-equivalent $kQ$-module $Y=X\oplus {^{\sigma}X}\oplus
\cdots\oplus {^{\sigma^{m-1}}X}$, there are induced $r=n/m$
indecomposable $kQ\#k\langle\sigma\rangle$-modules, which are
$(kQ\#k\langle\sigma\rangle)\otimes_{kQ\#k\langle\sigma^{m}\rangle}(L_{i}\otimes_{k}X)$,
$i=1,\cdots,r$;

$(II)$ For any indecomposable $kQ\#k\langle\sigma\rangle$-module
$Y$, there exists an indecomposable $kQ$-module $X$, so then apply
$(I)$, and there exists a unique $j\in \{1,2,\cdots, r\}$, such
that $Y\cong
(kQ\#k\langle\sigma\rangle)\otimes_{kQ\#k\langle\sigma^{m}\rangle}(L_{j}\otimes_{k}X)$.

We end the section by giving the relation between simple,
projective and injective modules between  mod$kQ$ and
mod$kQ\#k\langle\sigma\rangle$.

\begin{lemma}\label{lem 4.10} Let $H$ be a finite dimensional semisimple Hopf
  algebra and $A$ a finite dimensional $H$-module algebra. For a left
  $A\#H$-module $I$, if $I$ is an injective $A$-module, then $I$ is an
  injective $A\#H$-module.
  \end{lemma}

  \emph{Proof.} For an $A\#H$-module $M, N$, let $g:M\rightarrow N$
  and $h:M\rightarrow I$ be two $A\#H$-module homomorphisms such
  that $g$ is injective. In order to prove that $I$ is injective
  as an $A\#H$-module, it is enough to find an $\tilde{f}\in
  Hom_{A\#H}(N,I)$ satisfying $h=\tilde{f}g$. Since $I$ is
  injective as an $A$-module, there is an $f\in Hom_{A}(N,I)$ such
  that $h=fg$, where we consider $A\#H$-modules as $A$-modules in
  the natural way. Define $\tilde{f}(n)=\sum_{(t)}S(t_{1})\cdot f(t_{2}\cdot n)$
  for $n\in N$, where $t$ is a non-zero right integral with
  $\varepsilon(t)=1$. Then $\tilde{f}$ is $A\#H$-linear by
  Proposition 2 in \cite{CF}, and $h=\tilde{f}g$ since $\tilde{f}g(m)=\sum_{(t)}S(t_{1})\cdot f(t_{2}\cdot g(m))=
  \sum_{(t)}S(t_{1})\cdot f(g(t_{2}\cdot m))=\sum_{(t)}S(t_{1})\cdot fg(t_{2}\cdot m)=\sum_{(t)}S(t_{1})\cdot h(t_{2}\cdot m)
  =\sum_{(t)}S(t_{1})\cdot (t_{2}\cdot h(m))=(\sum_{(t)}S(t_{1})t_{2})\cdot
  h(m)=\varepsilon(t)h(m)=h(m)$.$\;\;\;\;\;\blacksquare$

Recall that $t\in H$ is a non-zero right integral, if
$th=\varepsilon(h)t,$ for any $h\in H$. Since $H$ is semisimple
Hopf algebra, there must exist a non-zero right integral such that
$\varepsilon(t)=1$. For details, see \cite{M}.

  \begin{theorem} Let $X$ be a $kQ\#k\langle\sigma\rangle$-module, then:
\begin{itemize}
  \item[$(i)$] $X$ is simple if and only if there exists a simple $kQ$-module
  $S$, such that $X$ is isomorphic to one of the $kQ\#k\langle\sigma\rangle$-modules
 induced from the indecomposable $\langle\sigma\rangle$-equivalent
 $kQ$-module $\oplus_{i=1}^{m-1}{^{\sigma^{i}}S}$.

  \item[$(ii)$] $X$ is  projective if and only if there exists an indecomposable projective $kQ$-module
  $P$, such that $X$ is isomorphic to one of the $kQ\#k\langle\sigma\rangle$-modules
  induced  from the indecomposable $\langle\sigma\rangle$-equivalent
 $kQ$-module
  $\oplus_{i=1}^{m-1}{^{\sigma^{i}}P}$.

  \item[$(iii)$] $X$ is injective if and only if there exists an indecomposable injective $kQ$-module
  $I$, such that $X$ is isomorphic to one of the $kQ\#k\langle\sigma\rangle$-modules
 induced  from the indecomposable $\langle\sigma\rangle$-equivalent
 $kQ$-module
  $\oplus_{i=1}^{m-1}{^{\sigma^{i}}I}$.
  \end{itemize}

  \end{theorem}

  \emph{Proof.} According to Theorem~\ref{thm 4.8} and Theorem~\ref{thm
  4.9}, we need only to prove that for a $kQ\#k\langle\sigma\rangle$-module, $X$ is a semisimple (projective, injective)
  $kQ\#k\langle\sigma\rangle$-module if and only if $X$ is a
  semisimple (projective, injective) $kQ$-module.

  $(i)$ Any $kQ\#k\langle\sigma\rangle$-module $X$
  induced
  from an indecomposable $\langle\sigma\rangle$-equivalent
  $kQ$-module
  $\bigoplus_{i=1}^{m-1}{^{\sigma^{i}}S}$ is a simple
  $kQ\#k\langle\sigma\rangle$-module since the dimension of vector space
  $X_{i}$  is 0 or 1,for any $i\in Q_{0}$. Additionally, any simple $kQ\#k\langle\sigma\rangle$-module $X$ is a semisimple
  $kQ$-module as in \cite{RR}.

  $(ii)$ By Lemma 3.1.7 in \cite{L}, for a $kQ\#k\langle\sigma\rangle$-module, $X$ is a projective
  $kQ\#k\langle\sigma\rangle$-module if and only if $X$ is a projective
  $kQ$-module.

 $(iii)$ By Lemma~\ref{lem 4.10}, for a $kQ\#k\langle\sigma\rangle$-module, $X$ is
 an
 injective
  $kQ\#k\langle\sigma\rangle$-module if and only if $X$ is an
  injective
  $kQ$-module.

  \section{Applications}

In this section, we apply the results we have gotten by giving
some examples. Let $\xi$ be the $n$-the primitive root of $1$.

\begin{example} Given a quiver $Q$, $\sigma\in Aut(Q)$ of order $n$. In this example, we are going to construct all
 simple $kQ\#k\langle\sigma\rangle$-modules in a concrete way. Let
 $S$ be a simple $kQ$-module with $m$ minimal such that $^{\sigma^{m}}S\cong
 S$ (in fact $^{\sigma^{m}}S=
 S$ since $S$ is simple), and let $r=n/m$.

 Let $\mathcal{S}^{(l)}, l\in \{0,1,\cdots, r-1\}$, as $kQ$-module, is an indecomposable
 $\langle\sigma\rangle$-equivalent $kQ$-module $S\oplus
 {^{\sigma}S}\oplus\cdots\oplus
 {^{\sigma^{m-1}}S}$, and $\sigma'$s action on $\mathcal{S}^{(l)}=S\oplus
 {^{\sigma}S}\oplus\cdots\oplus
 {^{\sigma^{m-1}}S}$ is defined by
\[\sigma(x_{0},x_{1},\cdots,x_{m-1})=(\xi^{ml}x_{m-1},x_{0},\cdots,x_{m-2}),x_{i}\in S.\]

Claim $1:$ $p(x)=\sigma^{-m}(p)(x),$ for any $x\in S$.

\emph{Proof.} For a simple $kQ$-module $S$, there exists a unique
$i\in Q_{0}$, such that for any $x\in S, q\in Q_{0}\cup
Q_{1}\setminus \{e_{i}\}$, $e_{i}(x)=x, q(x)=0$. Since
$\sigma^{m}S= S$, then $\sigma^{m}(e_{i})=e_{i}, \sigma^{q}\in
Q_{0}\cup Q_{1}\setminus \{e_{i}\},$ for any $q\in Q_{0}\cup
Q_{1}\setminus \{e_{i}\}$. $\blacksquare$

Claim $2:$  For any $l\in \{0,1,\cdots, r-1\},\mathcal{S}^{(l)}$
is a $kQ\#k\langle\sigma\rangle$-module.

 \emph{Proof.} We need
only to prove that two equations $\sigma^{n}=1$ and
$p\#\sigma=\sigma \sigma^{-1}(p)$ are satisfied as actions on
$\mathcal{S}^{(l)}$:
\begin{eqnarray*}\sigma^{n}(x_{0},x_{1},\cdots,x_{m-1})&=&\xi^{mr}(x_{0},x_{1},\cdots,
 x_{m-1})\\&=&(x_{0},x_{1},\cdots,x_{m-1}),\\
(p\#\sigma)(x_{0},x_{1},\cdots,x_{m-1})&=& p
(\xi^{ml}x_{m-1},x_{0},\cdots,x_{m-2})
\\&=&(\xi^{ml}p(x_{m-1}),\sigma^{-1}(p)(x_{0}),\cdots, \sigma^{-(m-1)}(p)(x_{m-2}))\\
&=& \sigma(\sigma^{-1}(p)(x_{0}), \sigma^{-2}(p)(x_{1}),\cdots, p(x_{m-1}))\\
 &\stackrel{Claim 1}=& \sigma(\sigma^{-1}(p)(x_{0}), \sigma^{-2}(p)(x_{1}),\cdots, \sigma^{-m}(p)(x_{m-1}))
\\&=& (\sigma \sigma^{-1}(p))(x_{0},x_{1},\cdots,x_{m-1}).\;\;\;\; \blacksquare
\end{eqnarray*}

Claim $3:$ If $l_{1}\neq l_{2}\in \{0,1,\cdots, r-1\},$ then
$\mathcal{S}^{(l_{1})}\ncong \mathcal{S}^{(l_{2})}$ as
$kQ\#k\langle\sigma\rangle$-modules.

\emph{Proof.} Simply  let $l_{1}=0, l_{2}=l,$ for some $l\in
\{1,2,\cdots, r-1\}$. Otherwise, there exists a
$kQ\#k\langle\sigma\rangle$-isomorphism
$F:\mathcal{S}^{(0)}\rightarrow \mathcal{S}^{(l)},$ denoted
$F(x,0,\cdots,0)=(F(x)_{0},F(x)_{1},\cdots,F(x)_{m-1})$.
\begin{eqnarray*}F(\sigma(0,\cdots,0,x))&=&F(x,0,\cdots, 0)\\&=&
(F(x)_{0},F(x)_{1},\cdots,F(x)_{m-1}),\\
\sigma(F(0,\cdots,0,x))&=&\sigma
(F(\sigma^{m-1}(x,0,\cdots,0)))
\\&=& \sigma^{m}F(x,0,\cdots,0)
\\&=&\sigma^{m}(F(x)_{0},F(x)_{1},\cdots,F(x)_{m-1})
\\&=&\xi^{ml}(F(x)_{0},F(x)_{1},\cdots,F(x)_{m-1}).
\end{eqnarray*}

From $F\sigma=\sigma F,$ we get $\xi^{ml}=1,$ which is
contradicted since $l\in\{1,2,\cdots,r-1\}$ and $\xi$ is a
primitive root of $1$. $\blacksquare$

So $\{\mathcal{S}^{(0)},\mathcal{S}^{(1)},\cdots,
\mathcal{S}^{(r-1)}\}$ are exactly $r$ non-isomorphic
$kQ\#k\langle\sigma\rangle$-modules induced from an indecomposable
$\langle\sigma\rangle$-equivalent $kQ$-module $S\oplus
{^{\sigma}S}\oplus \cdots\oplus{^{\sigma^{m-1}}S}$.

\end{example}

\begin{example}
Given a quiver $Q$, $\sigma\in Aut(Q)$ of order $n$. In this
example, we are going to construct all
 indecomposable projective $kQ\#k\langle\sigma\rangle$-modules in a concrete way. Let
 $P$ be an indecomposable $kQ$-module with $m$ minimal such that $^{\sigma^{m}}P\cong
 P$, and let $r=n/m$.

 Let $\mathcal{P}^{(l)}, l\in \{0,1,\cdots, r-1\}$, as $kQ$-module, is an indecomposable
 $\langle\sigma\rangle$-equivalent $kQ$-module $P\oplus
 {^{\sigma}P}\oplus\cdots\oplus
 {^{\sigma^{m-1}}P}$, and $\sigma'$s action on $\mathcal{P}^{(l)}=P\oplus {^{\sigma}P}\oplus \cdots\oplus
 {^{\sigma^{m-1}}P}$ is defined by
\[\sigma(x_{0},x_{1},\cdots,x_{m-1})=(\xi^{ml}\sigma^{m}(x_{m-1}),x_{0},\cdots,x_{m-2}),x_{i}\in P.\]

Claim $1:$ The action of $\sigma$ is well-defined due to
$\sigma^{m}(x)\in P$ for any $x\in P$.

\emph{Proof.} For an indecomposable projective $kQ$-module $P$,
there exists a unique $i\in Q_{0}$, such that $P=kQe_{i}$. Since
$kQ$ is a $k\langle\sigma\rangle$-module algebra, we have
$\sigma^{m}(P)=\sigma^{m}(kQ)\sigma^{m}(e_{i})=kQ\sigma^{m}(e_{i})$,
then we need only to prove $\sigma^{m}(e_{i})=e_{i}$. It is
clearly true by considering the simple module $S=P/rP$.
$\blacksquare$

Claim $2:$  For any $l\in \{0,1,\cdots, r-1\},\mathcal{P}^{(l)}$
is a $kQ\#k\langle\sigma\rangle$-module.

 \emph{Proof.} We need
only to prove that two equation $\sigma^{n}=1$ and
$p\#\sigma=\sigma \sigma^{-1}(p)$, are satisfied as actions on
$\mathcal{S}^{(l)}$:
\begin{eqnarray*}\sigma^{n}(x_{0},x_{1},\cdots,x_{m-1})&=&\xi^{mr}(\sigma^{mr}(x_{0}),\sigma^{mr}(x_{1}),\cdots,
 \sigma^{mr}(x_{m-1}))\\&=&(x_{0},x_{1},\cdots,x_{m-1}),\\
(p\#\sigma)(x_{0},x_{1},\cdots,x_{m-1})&=& p
(\xi^{ml}\sigma^{m}(x_{m-1}),x_{0},\cdots,x_{m-2})
\\&=&(\xi^{ml}p\sigma^{m}(x_{m-1}),\sigma^{-1}(p)(x_{0}),\cdots, \sigma^{-(m-1)}(p)(x_{m-2}))\\
&=& (\xi^{ml}\sigma^{m}\sigma^{-m}(p)(x_{m-1}), \sigma^{-1}(p)(x_{1}),\cdots, \sigma^{-(m-1)}(p)(x_{m-2}))\\
 &=& \sigma(\sigma^{-1}(p)(x_{0}), \sigma^{-2}(p)(x_{1}),\cdots, \sigma^{-m}(p)(x_{m-1}))
\\&=& (\sigma \sigma^{-1}(p))(x_{0},x_{1},\cdots,x_{m-1}).\;\;\;\; \blacksquare
\end{eqnarray*}

Claim $3:$ If $l_{1}\neq l_{2}\in \{0,1,\cdots, r-1\},$ then
$\mathcal{P}^{(l_{1})}\ncong \mathcal{P}^{(l_{2})}$ as
$kQ\#k\langle\sigma\rangle$-modules.

\emph{Proof.} Simply we let $l_{1}=0, l_{2}=l,$ for some $l\in
\{1,2,\cdots, r-1\}$. Otherwise, there exists a
$kQ\#k\langle\sigma\rangle$-isomorphism
$F:\mathcal{P}^{(0)}\rightarrow \mathcal{P}^{(l)},$ denoted
$F(x,0,\cdots,0)=(F(x)_{0},F(x)_{1},\cdots,F(x)_{m-1})$.
\begin{eqnarray*}F(\sigma(0,\cdots,0,x))&=&F(\sigma^{m}(x),0,\cdots, 0)\\&=&
(F(\sigma^{m}(x))_{0},F(\sigma^{m}(x))_{1},\cdots,F(\sigma^{m}(x))_{m-1}),\\
\sigma(F(0,\cdots,0,x))&=&\sigma (F(\sigma^{m-1}(x,0,\cdots,0)))
\\&=& \sigma^{m}F(x,0,\cdots,0)
\\&=&\sigma^{m}(F(x)_{0},F(x)_{1},\cdots,F(x)_{m-1})
\\&=&\xi^{ml}(\sigma^{m}(F(x)_{0}),\sigma^{m}(F(x)_{1}),\cdots,\sigma^{m}(F(x)_{m-1})).
\end{eqnarray*}

From $F\sigma=\sigma F,$ particularly, $F\sigma(e_{i})=\sigma
F(e_{i})$, we get $\xi^{ml}=1,$ which is contradicted since
$l\in\{1,2,\cdots,r-1\}$ and $\xi$ is a primitive root of $1$.
$\blacksquare$

So $\{\mathcal{P}^{(0)},\mathcal{P}^{(1)},\cdots,
\mathcal{P}^{(r-1)}\}$ are the $r$ non-isomorphic
$kQ\#k\langle\sigma\rangle$-modules induced from an indecomposable
$\langle\sigma\rangle$-equivalent $kQ$-module $P\oplus
{^{\sigma}P}\oplus\cdots\oplus {^{\sigma^{m-1}}P}$.

\end{example}

\begin{example}
Let $Q$ be the Kronecker quiver

\begin{center}
\begin{picture}(200,20)(0,0)
 \put(50,10){\makebox(0,0){$ a_{0}$}}
\put(75,18){\makebox(0,0){$ \alpha_{1}$}}
\put(75,-2){\makebox(0,0){$\alpha_{0}$}}
 \put(90,15){\vector(-1,0){30}}
 \put(90,5){\vector(-1,0){30}}
\put(100,10){\makebox(0,0){$a_{1}$}}
\end{picture}
\end{center}

Let $p(l),i(l),l\in \mathbf{N},r_{\lambda}(l),r_{\infty}(l),l\in
\mathbf{N}\setminus \{0\},\lambda\in k$ be the $kQ$-modules
defined by

\begin{center}
\begin{picture}(200,20)(80,20)
\put(20,10){\makebox(0,0){$p(l):$}}
 \put(50,10){\makebox(0,0){$ k^{l+1}$}}
\put(75,25){\makebox(0,0){$ [\;_{0}^{\mathbf{I}_{l}}]$}}
\put(75,-5){\makebox(0,0){$[\;_{\mathbf{I}_{l}}^{0}]$}}
 \put(90,15){\vector(-1,0){30}}
 \put(90,5){\vector(-1,0){30}}
\put(100,10){\makebox(0,0){$k^{l}$}}
\put(220,10){\makebox(0,0){$i(l):$}}
 \put(250,10){\makebox(0,0){$ k^{l}$}}
\put(275,25){\makebox(0,0){$ [\mathbf{I}_{l},0]$}}
\put(275,-5){\makebox(0,0){$[0,\mathbf{I}_{l}]$}}
 \put(290,15){\vector(-1,0){30}}
 \put(290,5){\vector(-1,0){30}}
\put(300,10){\makebox(0,0){$k^{l+1}$}}
\put(20,-60){\makebox(0,0){$r_{\lambda}(l):$}}
 \put(50,-60){\makebox(0,0){$ k^{l}$}}
\put(75,-45){\makebox(0,0){$ J_{\lambda}(l)$}}
\put(75,-75){\makebox(0,0){$\mathbf{I}_{l}$}}
 \put(90,-55){\vector(-1,0){30}}
 \put(90,-65){\vector(-1,0){30}}
\put(100,-60){\makebox(0,0){$k^{l}$}}
\put(220,-60){\makebox(0,0){$r_{\infty}(l):$}}
 \put(250,-60){\makebox(0,0){$ k^{l}$}}
\put(275,-45){\makebox(0,0){$\mathbf{I}_{l}$}}
\put(275,-75){\makebox(0,0){$J_{0}(l)$}}
 \put(290,-55){\vector(-1,0){30}}
 \put(290,-65){\vector(-1,0){30}}
\put(300,-60){\makebox(0,0){$k^{l}$}}
\end{picture}\\
\end{center}
${}$\\

${}$\\

${}$\\

where $J_{\lambda}(l),\lambda\in k$ is the $l\times l$ Jordan
block

\[\left(%
\begin{array}{cccc}
  \lambda & 1 &  &  \\
   & \ddots & \ddots &  \\
   &  & \lambda & 1 \\
   &  &  & \lambda \\
\end{array}%
\right)\]

It is known in  \cite{GR, KL} that $\{p(l),i(l),l\in
\mathbf{N},r_{\lambda}(l),r_{\infty}(l),l\in \mathbf{N}\setminus
\{0\},\lambda\in k\}$ classify all indecomposable $kQ$-modules up
to isomorphism.

Let $P(l)^{(0)}, l\in \mathbf{N}$, be the
$kQ\#k\langle\sigma\rangle$-module that is $p(l)$ as $kQ$-module
and $\sigma'$s action is defined by

\[(\left(%
\begin{array}{cccc}
   &  &  & 1 \\
   &  & \cdots &  \\
   & 1 &  &  \\
  1 &  &  &  \\
\end{array}%
\right)_{l+1,l+1},\left(%
\begin{array}{cccc}
   &  &  & 1 \\
   &  & \cdots &  \\
   & 1 &  &  \\
  1 &  &  &  \\
\end{array}%
\right)_{l,l}).\]

Let $P(l)^{(1)}, l\in \mathbf{N}$, be the
$kQ\#k\langle\sigma\rangle$-module that is $p(l)$ as $kQ$-module
and $\sigma'$s action is defined by

\[(-\left(%
\begin{array}{cccc}
   &  &  & 1 \\
   &  & \cdots &  \\
   & 1 &  &  \\
  1 &  &  &  \\
\end{array}%
\right)_{l+1,l+1},-\left(%
\begin{array}{cccc}
   &  &  & 1 \\
   &  & \cdots &  \\
   & 1 &  &  \\
  1 &  &  &  \\
\end{array}%
\right)_{l,l}).\]

Let $I(l)^{(0)}, l\in \mathbf{N}$, be the
$kQ\#k\langle\sigma\rangle$-module that is $i(l)$ as $kQ$-module
and $\sigma'$s action is defined by

\[(\left(%
\begin{array}{cccc}
   &  &  & 1 \\
   &  & \cdots &  \\
   & 1 &  &  \\
  1 &  &  &  \\
\end{array}%
\right)_{l,l},\left(%
\begin{array}{cccc}
   &  &  & 1 \\
   &  & \cdots &  \\
   & 1 &  &  \\
  1 &  &  &  \\
\end{array}%
\right)_{l+1,l+1}).\]

Let $I(l)^{(1)}, l\in \mathbf{N}$, be the
$kQ\#k\langle\sigma\rangle$-module that is $i(l)$ as $kQ$-module
and $\sigma'$s action is defined by

\[(-\left(%
\begin{array}{cccc}
   &  &  & 1 \\
   &  & \cdots &  \\
   & 1 &  &  \\
  1 &  &  &  \\
\end{array}%
\right)_{l,l},-\left(%
\begin{array}{cccc}
   &  &  & 1 \\
   &  & \cdots &  \\
   & 1 &  &  \\
  1 &  &  &  \\
\end{array}%
\right)_{l+1,l+1}).\]

Let $R_{(0,\infty)}(l),l\in \mathbf{N}\setminus \{0\}$, be the
$kQ\#k\langle\sigma\rangle$-module that is  $r_{0}(l)\oplus
r_{\infty}(l)$ as $kQ$-modules and $\sigma'$s action is defined by
\[ \sigma((x_{0},x_{1}),(y_{0},y_{1}))=((y_{0},y_{1}),(x_{0},x_{1})),\forall x_{0},x_{1},y_{0},y_{1}\in k^{l}.
\]

Let $R_{(\lambda,\lambda^{-1})}(l),\lambda\in k\setminus \{0\},
l\in \mathbf{N}\setminus \{0\}$, be the
$kQ\#k\langle\sigma\rangle$-module that is $r_{\lambda}(l)\oplus
r_{\lambda^{-1}}(l)$ as $kQ$-modules and $\sigma'$s action is
defined by
\[\sigma((x_{0},x_{1}),(y_{0},y_{1}))=((B_{l}^{-1}(y_{0}), A_{l}^{-1}(y_{1})),(B_{l}(x_{0}),A_{l}(x_{1}))),\forall x_{0},x_{1},y_{0},y_{1}\in k^{l},\]
where $B_{l}=(b_{ij})_{l\times l},A_{l}=(a_{ij})_{l\times l}\in
M_{l\times l}(k)$ satisfy
\begin{eqnarray*} b_{ij}&=&0, \;\; i>j,\\
b_{il}&=& 0,\;\; 1\leq i<l,\\
b_{ll}&=& 1, \\
b_{ij}&=& -(b_{i-1,j}\lambda+b_{i-1,j+1}\lambda^{2}),\;\;
1\leq i<l, i\leq j<l.\\
a_{ij}&=&0, \;\; i>j,\\
a_{ll}&=& \lambda, \\
a_{ij}&=& -(a_{i-1,j}\lambda+a_{i-1,j+1}\lambda^{2}),\;\; 1\leq
i<l, i\leq j<l,\\
a_{il}&=& -a_{i-1,l}\lambda,\;\; 1\leq i<l.
\end{eqnarray*}

$\;$

It is easy to see that
\[\{P(l)^{(0)},P(l)^{(1)},I(l)^{(0)},I(l)^{(1)},l\in
\mathbf{N},R_{(0,\infty)}(l),R_{(\lambda,\lambda^{-1})}(l),l\in
\mathbf{N}\setminus \{0\},\lambda \in k\setminus \{0\} \}\]
classify all indecomposable $kQ\#k\langle\sigma\rangle$-modules up
to isomorphism.
\end{example}

\emph{\bf Acknowledgement.}  The paper was ultimately  finished
during the first author's visit at Kansas State University. The
first author thanks The China Council Scholarship  for the
financial support of her
 visit to Kansas State University. The authors thank Professor
Zongzhu Lin of Kansas State University for his discussions and
suggestions.

\end{document}